\documentclass[a4paper]{amsart}
\usepackage{amssymb}
\usepackage[all]{xy}

\newcommand{\fun}[3]{#1\colon #2\rightarrow #3}

\newcommand{\set}[2]{\{#1:\, #2\}}

\newcommand{\psr}[2]{#1[[#2]]}

\newcommand{\abs}[1]{\lvert#1\rvert}

\newcommand{\cc}[1]{\mathcal{#1}}

\newcommand{\dotsk}[0]{,\dots,}

\newcommand{\bbN}[0]{\mathbb{N}}
\newcommand{\bbZ}[0]{\mathbb{Z}}
\newcommand{\bbQ}[0]{\mathbb{Q}}

\newcommand{\ccP}[0]{\mathcal{P}}

\newcommand{\per}[0]{\Sigma}                          
\newcommand{\schur}[0]{Schur}                         
\newcommand{\ifc}[1]{C_{#1}}

\DeclareMathOperator{\Disjunion}{\overset{\centerdot}{\bigcup}}

\DeclareMathOperator{\ext}{\bigwedge}                             

\DeclareMathOperator{\rep}{R}                          
\DeclareMathOperator{\res}{res}                         
\DeclareMathOperator{\ind}{ind}                         
\DeclareMathOperator{\sgn}{sgn}
\DeclareMathOperator{\bur}{\mathcal{B}}                 
\DeclareMathOperator{\btr}{h}

\numberwithin{equation}{section}
\newtheorem{theorem}[equation]{Theorem}
\newtheorem{lemma}[equation]{Lemma}

\newtheorem{proposition}[equation]{Proposition}
\newtheorem{definition}[equation]{Definition}

\theoremstyle{remark}
\newtheorem*{remark}{Remark}
\newtheorem*{acknowledgements}{Acknowledgment}

\begin{document}

\title[A note on the $\lambda$-structure on the Burnside ring]{A note on the $\lambda$-structure on the Burnside ring}
\author{Karl R\"okaeus}
\address{Karl R\"okaeus \\ Matematiska institutionen \\ Stockholms universitet \\ SE-106 91 Stockholm \\ Sweden}
\email{karlr@math.su.se}
\date{August 10, 2007}

\begin{abstract}
Let $G$ be a finite group and let $S$ be a $G$-set. The Burnside ring of $G$ has a natural structure of a $\lambda$-ring, $\{\lambda^n\}_{n\in\bbN}$. However, a priori $\lambda^n(S)$ can only be computed recursively, by first computing $\lambda^1(S),\dots,\lambda^{n-1}(S)$. In this paper we establish an explicit formula, expressing $\lambda^n(S)$ as a linear combination of classes of $G$-sets.
\end{abstract}

\maketitle

\section{Introduction}
We use $\bur(G)$ to denote the Burnside ring of the finite group $G$. Recall that, as an abelian group, $\bur(G)$ is free on $\{[S]\}_{S\in R}$, where $R$ is a set of representatives of the isomorphism classes of transitive $G$-sets, and that its rank equals the number of conjugacy classes of subgroups of $G$. When $f$ is a function on $\bur(G)$, we write $f(S)$ for $f([S])$.

There is a $\lambda$-structure on $\bur(G)$, $\{\lambda^n\}_{n\in\bbN}$, defined as the opposite structure of $\{\sigma^n\}_{n\in\bbN}$, where $\sigma^n(S)$ is the class of the $n$th symmetric power of $S$.\footnote{Let $\{\sigma^n\}_{n\in\bbN}$ be a $\lambda$-structure on the ring $R$ and define $\sigma_t(x):=\sum_{i\geq0}\sigma^i(x)t^i\in\psr{R}{t}$. The $\lambda$-structure opposite to $\{\sigma^n\}$ is defined by $\sigma_t(x)\cdot\lambda_{-t}(x)=1\in\psr{R}{t}$, where $\lambda_t(x):=\sum_{i\geq0}\lambda^i(x)t^i$.} It should be considered the natural $\lambda$-structure on $\bur(G)$, one reason for this being that there is a canonical homomorphism to the ring of rational representations of $G$, $\btr\colon\bur(G)\to\rep_{\bbQ}(G)$, defined by $\btr(S)=\bigl[\bbQ[S]\bigr]$, and the given $\lambda$-structure on $\bur(G)$ makes $\btr$ into a $\lambda$-homomorphism. (Note however that this $\lambda$-structure is non-special.)

The implicit nature of the definition of the $\lambda$-structure on $\bur(G)$ makes it hard to compute with. The main result of this paper is a closed formula for $\lambda^i(S)$, where $S$ is any $G$-set. To state it we first introduce some notation:
\begin{definition}
Let $\mu\vdash i$, i.e., $\mu$ is a partition of $i$. We use $\ell(\mu)$ to denote the length of $\mu$. Also, if $\mu=(\mu_1\dotsk \mu_l)$, where $\mu_1=\dots=\mu_{\alpha_1}>\mu_{\alpha_1+1}=\dots=\mu_{\alpha_1+\alpha_2}>\dots>\mu_{l-\alpha_{l'}+1}=\dots=\mu_l$ we define the tuple $\alpha(\mu):=(\alpha_1\dotsk\alpha_{l'})$, and write $\tbinom{l(\mu)}{\alpha(\mu)}$ for $\tfrac{l!}{\alpha_1!\cdots\alpha_{l'}!}$. Finally, for $S$ a $G$-set of cardinality $\geq i$ we define $\mathcal{P}_\mu(S)$ to be the $G$-set consisting of $\ell(\mu)$-tuples of pairwise disjoint subsets of $S$, where the first one has cardinality $\mu_1$, and so on.
\end{definition}
Using this notation we can express $\lambda^n(S)$, for any $G$-set $S$, as a linear combination of classes of $G$-sets:
\begin{equation}\label{6}
\lambda^i(S)=(-1)^i\sum_{\mu\vdash i}(-1)^{\ell(\mu)}\tbinom{\ell(\mu)}{\alpha(\mu)}[\mathcal{P}_\mu(S)]\in\bur(G)
\end{equation}
when $i\leq\abs{S}$, and $\lambda^i(S)=0$ when $i>\abs{S}$. When $S$ is transitive, so are the $\ccP_{\mu}(S)$, so in this case, \eqref{6} expresses $\lambda^i(S)$ as a linear combination of different elements from the standard basis of $\bur(G)$.

We will prove \eqref{6} in the following way: First, recall that a group homomorphism $\phi\colon H\to G$ gives rise to a $\lambda$-homomorphism $\res_H^G\colon\bur(G)\to\bur(H)$ by restricting the action on a $G$-set $S$ to an $H$-action via $\phi$. Now let $G$ be a finite group and let $S$ be a $G$-set of cardinality $n$. By choosing an enumeration of $S$ we get a homomorphism $G\to\per_n$, the symmetric group on $\{1\dotsk n\}$. Let $\res_G^{\per_n}$ be the corresponding restriction homomorphism (which is independent of the chosen enumeration). We have that $\res_G^{\per_n}\bigl(\{1\dotsk n\}\bigr)=[S]$, hence $\res_G^{\per_n}\bigl(\lambda^i(\{1\dotsk n\})\bigr)=\lambda^i(S)$. Also, writing $\ccP^{(n)}_\mu$ for $\ccP_\mu(\{1\dotsk n\})$, we see that $\res_G^{\per_n}(\mathcal{P}^{(n)}_\mu)=[\mathcal{P}_\mu(S)]$. Hence, to prove \eqref{6} it suffices to prove it in the special case when $G=\per_n$ and $S=\{1\dotsk n\}$, i.e.,
\begin{equation}\label{3}
\lambda^i\bigl(\{1\dotsk n\}\bigr)=(-1)^i\sum_{\mu\vdash i}(-1)^{\ell(\mu)}\tbinom{\ell(\mu)}{\alpha(\mu)}[\mathcal{P}^{(n)}_\mu]\in\bur(\per_n).
\end{equation}
The validity of \eqref{3} will be established in Theorem \ref{80}.

We have not been able to derive \eqref{3} completely inside of the Burnside ring. Instead we have to use the canonical $\lambda$-homomorphism $\btr\colon\bur(\per_n)\to\rep_{\bbQ}(\per_n)$ to move some of the computations to the rational representation ring, whose $\lambda$-structure is much easier to work with. However, there is a problem in that $\btr$ is not injective for $\per_n$. In Section \ref{1} we will therefore introduce a subring $\schur_n$ of $\bur(\per_n)$, with the property that the restriction of $\btr\colon\bur(\per_n)\to\rep_{\bbQ}(\per_n)$ to $\schur_n$ is injective.

In Section \ref{2} we then establish \eqref{3}. The technique of passing to the representation ring will be used at a crucial place, to prove Lemma \ref{73}.

In a forthcoming paper, \cite{rokaeusGro}, we give an application of \eqref{3}: We use it to derive a formula for the class of a certain torus in the Grothendieck ring of varieties.

An introduction to $\lambda$-rings, representation rings and the Burnside ring is given in \cite{MR0364425}. The standard reference for $\lambda$-rings is \cite{MR0244387}, whereas representation rings are extensively studied in \cite{MR0450380}.
\begin{acknowledgements}
The author is grateful to Professor Torsten Ekedahl for valuable discussions and comments on the manuscript.
\end{acknowledgements}

\section{The Schur subring of $\bur(\per_n)$}\label{1}
Recall that we write $\ccP_\mu^{(n)}$ for $\ccP_\mu(\{1\dotsk n\})\in\bur(\per_n)$.

Let $S$ be a $\per_n$-set. We say that $S$ is a \emph{Schur set} if for every $s\in S$, $(\per_n)_s$, the stabilizer subgroup of $s$, is a Schur subgroup, i.e., it is the stabilizer of some partition of $\{1\dotsk n\}$. Equivalently, any transitive component of $S$ is isomorphic to $\ccP_\mu^{(n)}$ for some $\mu\vdash n$.
\begin{definition}
$\schur_n$ is the subgroup of $\bur(\per_n)$ generated by the Schur sets. 
\end{definition}
Equivalently, this means that $\schur_n\subset\bur(\per_n)$ is the free subgroup on $\{[\ccP_\mu]\}_{\mu\vdash n}$. The reason for us to introduce $\schur_n$ is the next theorem:
\begin{theorem}
Let $\fun{\btr}{\bur(\per_n)}{\rep_{\bbQ}(\per_n)}$ be the canonical $\lambda$-ring homomorphism. The restriction of $\btr$ to $\schur_n$ is injective.
\end{theorem}
Even though this is a simple consequence of the injectivity of the character homomorphism from $\rep_{\bbQ}(\per_n)$ to the ring of symmetric polynomials, we have chosen to give a more direct proof:
\begin{proof}
For every $\mu\vdash n$, let $\sigma_\mu\in\per_n$ be an element in the conjugacy class determined by $\mu$ and let $\ifc{\sigma_{\mu}}\colon\rep_{\bbQ}(\per_n)\rightarrow\bbZ$ be the homomorphism defined by $V\mapsto\chi_V(\sigma_\mu)$, where $\chi_V$ is the character of $V$. This definition is independent of the choice of $\sigma_\mu$. Together the $\ifc{\sigma_{\mu}}$ give a homomorphism
$$\rep_{\bbQ}(\per_n)\rightarrow\prod_{\mu\vdash n}\bbZ,$$
and it suffices to show that the composition of this with the restriction of $\btr$ to $\schur_n$ is injective, i.e., that
\begin{align*}
\varphi\colon\schur_n&\rightarrow\prod_{\mu\vdash n}\bbZ\\
[T]&\mapsto\bigl(\abs{T^{\sigma_\mu}}\bigr)_{\mu\vdash n}
\end{align*}
is injective, where $T^{\sigma_\mu}$ is the set of points in $T$ fixed by $\sigma_\mu$. To do this, define a total ordering on the set of partitions of $n$ by $\mu>\mu'$ if $\mu_1=\mu'_1\dotsk\mu_{j-1}=\mu'_{j-1}$ and $\mu_j>\mu'_j$ for some $j$ (i.e., lexicographic order). We claim that $\abs{\mathcal{P}_{\mu}^{\sigma_\mu}}\neq0$, whereas $\mathcal{P}_{\mu'}^{\sigma_\mu}=\emptyset$ if $\mu>\mu'$. (Here and in the rest of this proof we write $\ccP_\mu$ for $\ccP_\mu^{(n)}$.)

For the first assertion, choose for example
$$\sigma_\mu=(1\dotsk\mu_1)(\mu_1+1\dotsk\mu_1+\mu_2)\cdots(n-\mu_{\ell(\mu)}+1\dotsk n).$$
Then
$$\bigl(\{1\dotsk\mu_1\},\{\mu_1+1\dotsk\mu_1+\mu_2\}\dotsk\{n-\mu_{\ell(\mu)}+1\dotsk n\}\big)\in\mathcal{P}_\mu$$
is fixed by $\sigma_\mu$.

For the second assertion, suppose $\mu'<\mu$ and $t=(T_1\dotsk T_l)\in\mathcal{P}_{\mu'}$, where $l=\ell(\mu')$. Suppose moreover that $t$ is fixed by $\sigma_\mu$. If now $\mu_1>\mu_2>\dots>\mu_{\ell(\mu)}$, then, with the same $\sigma_\mu$ as above, we must have $T_1=\{1\dotsk\mu_1\}\dotsk T_l=\{n-\mu_l+1\dotsk n\}$. (This is because $\mu_j\geq\mu'_j$ for every $j$ and if $1$ lies in $T_j$ then so does $\sigma_\mu(1)=2$, hence also $3\dotsk\mu_1$. So $T_j$ has cardinality at least $\mu_1$ and the only $\mu'_j$ that can be that big is $\mu'_1$.) But if $\mu$ and $\mu'$ differ in position $j$ it is impossible for $T_j$ to fulfill this since it has cardinality $\mu'_j<\mu_j$. In the general case, when we may have $\mu_j=\mu_{j+1}$, the above argument works the same only that we for example can have $T_1=\{\mu_1+1\dotsk\mu_1+\mu_2\}$ and $T_2=\{1\dotsk\mu_1\}$ if $\mu_1=\mu_2$.

We are now ready to prove that $\varphi$ is injective. Let $x=\sum_{\mu\vdash n}a_\mu[\mathcal{P}_\mu]$, where $a_\mu\in\bbZ$, and suppose that $x\neq0$. Choose the maximal $\mu_0$ such that $a_{\mu_0}\neq0$. Let $\varphi_{\mu_0}$ be the $\mu_0$th component of $\varphi$. Then
$$\varphi_{\mu_0}(x)
=\sum_{\mu\vdash n}a_\mu\abs{\mathcal{P}_\mu^{\sigma_{\mu_0}}}
=a_{\mu_0}\abs{\mathcal{P}_{\mu_0}^{\sigma_{\mu_0}}}\neq0;$$
hence $\varphi(x)\neq0$.
\end{proof}

We conclude this section by showing that $\schur_n$ is a subring of $\bur(\per_n)$.
\begin{proposition}\label{7}
$\schur_n$ is closed under multiplication.
\end{proposition}
\begin{proof}
We want to see what happens when we multiply $[S]$ and $[T]$, where $S$ and $T$ are Schur sets. Let $s\in S$ and $t\in T$. Then the stabilizers of $s$ and $t$ equal the stabilizers of partitions of $\{1\dotsk n\}$, which we denote $(S_1\dotsk S_k)$ and $(T_1\dotsk T_l)$, respectively. Then $\sigma\in\per_n$ is in the stabilizer of $(s,t)$ precisely when $\sigma S_i=S_i$ and $\sigma T_j=T_j$ for each $i,j$. Equivalently, $\sigma$ must preserve $S_i\cap T_j$ for each $i,j$. Hence $(\per_n)_{(s,t)}$ equals the stabilizer of the partition $\{S_i\cap T_j\}_{i,j}$. Consequently, it is a Schur subgroup, hence $S\times T$ is a Schur set. Therefore $\schur_n$ is closed under multiplication.
\end{proof}
\begin{remark}
$\schur_n$ is in general not a $\lambda$-ring since it is not closed under the $\lambda$-operations. For example, note that when $S$ is a $\per_n$-set we can represent the symmetric square of $S$ as the set of $2$-subsets of $S$. Now, let $S$ be the $\per_4$-set $\{1,2,3,4\}$ and consider $x:=\sigma^2\bigl(\sigma^2(S)\bigr)\in\bur(\per_4)$. An element of the underlying $\per_4$-set is $\{\{1,2\},\{3,4\}\}$ and the stabilizer $G$ of this element is generated by $\{(12),(34),(13)(24),(14)(23)\}$. The only partition that is stabilized by $G$ is the trivial one, so since $G$ does not equal $\per_4$ it fails to be the stabilizer of a partition. Hence $x\notin\schur_4$. Since $\lambda^2\bigl(\sigma^2(S)\bigr)=\bigl(\sigma^2(S)\bigr)^2-x$ it follows that this is not contained in $\schur_4$ either. But $\sigma^2(S)$ is in $\schur_4$, which is therefore not a $\lambda$-ring. This also gives an example showing that $\bur(\per_4)$ is not special. For if it were, then $y:=\lambda^2\bigl(\lambda^2(S)\bigr)$ would be a polynomial in $\lambda^i(S)$ for $i=1,2,3,4$, which all lie in $\schur_4$, so $y$ would also lie in $\schur_4$. But using the above one shows that $y\notin\schur_4$.
\end{remark}

\section{The $\lambda$-operations on $\bur(\per_n)$}\label{2}
We are now ready to start the investigation of how $\lambda^i$ acts on $\{1\dotsk n\}$, the goal being to obtain a closed formula for it. We will need some more definitions. We have only defined $\ccP_\mu^{(n)}$ when $\mu$ is a partition of $i\leq n$. More generally:
\begin{definition}
If $\alpha=(i_1\dotsk i_l)$ is any tuple of positive integers summing up to $i\leq n$ we define $\ccP_{\alpha}^{(n)}$ to be the $\per_n$-set of $l$-tuples of disjoint subsets of $\{1\dotsk n\}$, the first one having $i_1$ elements, and so on.
\end{definition}
We have $[\ccP^{(n)}_{\alpha}]=[\ccP^{(n)}_\mu]$, where $\mu$ is the $i$-tuple corresponding to $\alpha$. Also note that $[\ccP^{(n)}_{\alpha}]=[\ccP^{(n)}_{\alpha,n-i}]$, where we use $\ccP^{(n)}_{\alpha,n-i}$ to denote $\ccP^{(n)}_{(i_i\dotsk i_l,n-i)}$. Similarly, if $\beta=(j_1\dotsk j_k)$ we will write $\ccP^{(n)}_{\alpha,\beta}$ for $\ccP^{(n)}_{(i_1\dotsk i_l,j_1\dotsk j_k)}$.

Throughout this section we use the following notation:
\begin{align*}
s^{(n)}_i:=&\sigma^i(\{1\dotsk n\})\\
\ell^{(n)}_i:=&\lambda^i(\{1\dotsk n\})\in\bur(\per_n)
\end{align*}
We begin by giving a formula for $s^{(n)}_i$ which shows that it lies in $\schur_n$, and we then deduce from this that also $\ell^{(n)}_i$ is in $\schur_n$. Recall from the introduction that if $\mu=(\mu_1,\dotsk \mu_l)$, where $\mu_1=\dots=\mu_{\alpha_1}>\mu_{\alpha_1+1}=\dots=\mu_{\alpha_1+\alpha_2}>\dots>\mu_{j-\alpha_{l'}+1}=\dots=\mu_j$, then we define $\alpha(\mu):=(\alpha_1\dotsk\alpha_{l'})$.
\begin{proposition}\label{8}
We have
$$s_i^{(n)}=\sum_{\substack{\mu\vdash i:\\ \ell(\mu)\leq n}}[\cc{P}^{(n)}_{\alpha(\mu)}].$$
In particular, $s_i^{(n)}$ and $\ell_i^{(n)}$ are in $\schur_n$ for every $i$.
\end{proposition}
\begin{proof}
Identify $\{1\dotsk n\}$ with $\{x_1\dotsk x_n\}$. Then the symmetric $i$th power of $\{1\dotsk n\}$, the $\per_n$-set $\{1\dotsk n\}^i/\per_i$, is identified with the set of monomials
$$\set{x_1^{e_1}\cdots x_n^{e_n}}{e_1+\dots+ e_n=i}=\Disjunion_{\substack{e_1+\dots+e_n=i\\ e_1\geq e_2\dots\geq e_n\geq0}}\per_n\cdot x_1^{e_1}\cdots x_n^{e_n},$$
where the index set on the disjoint union can be identified with the set of $\mu\vdash i$ such that $\ell(\mu)\leq n$. Now let $e_1=\dots=e_{\alpha_1}>e_{\alpha_1+1}=\dots=e_{\alpha_1+\alpha_2}>\dots>e_{n-\alpha_l+1}=\dots=e_n$. Then
\begin{align*}
\per_n\cdot x_1^{e_1}\cdots x_n^{e_n}=&\per_n\cdot (x_1\cdots x_{\alpha_1})^{e_1}(x_{\alpha_1+1}\cdots x_{\alpha_1+\alpha_2})^{e_{\alpha_1+1}}\cdots (x_{n-\alpha_l+1}\cdots x_n)^{e_{n-\alpha_l+1}}\\
\simeq&\per_n\bigl(\{x_1\dotsk x_{\alpha_1}\},\{x_{\alpha_1+1}\dotsk x_{\alpha_1+\alpha_2}\}\dotsk\{x_{n-\alpha_l+1}\dotsk x_n\}\bigr)\\
\simeq&\cc{P}^{(n)}_{(\alpha_1\dotsk\alpha_l)}
\end{align*}
so the first part of the proposition follows.

To show that also $\ell_i^{(n)}\in\schur_n$ we use that, by definition,
\begin{equation*}
-(-1)^i\ell_i^{(n)}=\sum_{j=0}^{i-1}(-1)^j\ell_j^{(n)}s_{i-j}^{(n)}.
\end{equation*}
Since we know that $\schur_n$ is a ring, and that all $s_j^{(n)}$ and $\ell_1^{(n)}=[\cc{P}_1^{(n)}]$ are in $\schur_n$, it follows by induction that $\ell_i^{(n)}\in\schur_n$.
\end{proof}

Recall the definition of the induction map, analogous to that for representation rings: If $\phi\colon H\to G$ is a homomorphism of groups and $S$ is a $H$-set, then we can associate to it the $G$-set $G\times_HS$, i.e., the quotient of $G\times S$ by the equivalence relation $(g\cdot\phi(h),s)\sim(g,hs)$ for $(g,s)\in G\times S$ and $h\in H$, with a $G$-action given by $g'\cdot(g,s):=(g'g,s)$. This gives rise to the induction map $\ind_H^G\colon\bur(H)\rightarrow\bur(G)$, which is additive but not multiplicative. We will only use it in the case when $H$ is a subgroup of $G$. In this case, note that if we choose a set of coset representatives of $G/H$, say $R=\{g_1\dotsk g_r\}$, then we can represent $G\times_HS$ as $R\times S$ with $G$-action given by $g\cdot(g_i,s)=(g_j,hs)$, where $gg_i=g_jh$ for $h\in H$. The map $\btr\colon\bur(G)\to\rep_{\bbQ}(G)$ commutes with the induction maps if $H$ is a subgroup of $G$.

In the following two lemmas we show that $\ell_i^{(n)}$ and $[\ccP_\mu^{(n)}]$, where $\mu\vdash i$, are determined by $\ell_i^{(i)}$ and $[\ccP_\mu^{(i)}]$ respectively. For this we use the map $$\ind_{\per_i\times\per_{n-i}}^{\per_n}\circ\res_{\per_i\times\per_{n-i}}^{\per_i}\colon\bur(\per_i)\to\bur(\per_n)$$
which is constructed in the following way: We view $\per_i$ as the symmetric group on $\{1\dotsk i\}$ and embed it in $\per_n$, the symmetric group on $\{1\dotsk n\}$. Moreover we view $\per_{n-i}$ as the symmetric group on $\{i+1\dotsk n\}$. We then restrict from $\bur(\per_i)$ to $\bur(\per_i\times\per_{n-i})$ with respect to the projection $\per_i\times\per_{n-i}\to\per_i$ and we induce from $\bur(\per_i\times\per_{n-i})$ to $\bur(\per_n)$ with respect to the inclusion $(\tau,\rho)\mapsto\tau\rho=\rho\tau\colon\per_i\times\per_{n-i}\to\per_n$.
\begin{lemma}\label{74}
Let $\mu\vdash i$. For $n\geq i$, $\ind_{\per_i\times\per_{n-i}}^{\per_n}\circ\res_{\per_i\times\per_{n-i}}^{\per_i}\bigl(\ccP_\mu^{(i)}\bigr)=[\ccP_\mu^{(n)}]\in\bur(\per_n)$.
\end{lemma}
\begin{proof}
Let $R=\{\sigma_1\dotsk\sigma_r\}$, where $r=\tbinom{n}{i}$, be a system of coset representatives for $\per_n/(\per_i\times\per_{n-i})$. We know that $\per_n\times_{\per_i\times\per_{n-i}}\ccP_\mu^{(i)}$ can be identified with the set of pairs $(\sigma_j,t)$, where $\sigma_j\in R$ and $t=(T_1\dotsk T_l)\in\ccP_\mu^{(i)}$. From this set we define a map to $\ccP_\mu^{(n)}$ by
$$(\sigma_j,t)\mapsto\bigl(\sigma_jT_1\dotsk\sigma_jT_l,\sigma_j\{i+1\dotsk n\}\bigr).$$
This map is surjective for given $t'=(T'_1\dotsk T'_l,T'_{l+1})\in\ccP_{\mu}^{(n)}$, there is a $\sigma\in\per_n$ such that $\sigma\{1\dotsk\mu_1\}=T'_1\dotsk\sigma\{i-\mu_l+1\dotsk i\}=T'_l$. Let $\sigma_j\in R$ be such that $\sigma=\sigma_j\tau\rho$ where $(\tau,\rho)\in\per_i\times\per_{n-i}$. Then
$$\bigl(\sigma_j,\tau(\{1\dotsk\mu_1\}\dotsk\{i-\mu_l+1\dotsk i\})\bigr)\mapsto t'.$$
Since both sets have $n!/(\mu_1!\cdots\mu_l!(n-i)!)$ elements this is a bijection. Finally, the map is $G$-equivariant, hence it is an isomorphism.
\end{proof}

It is the following lemma that forces us to pass to the representation ring, for we have not been able to prove it directly in the Burnside ring.
\begin{lemma}\label{73}
Given $i$. For $n\geq i$ we have $\ind_{\per_i\times\per_{n-i}}^{\per_n}\circ\res_{\per_i\times\per_{n-i}}^{\per_i}\bigl(\ell_i^{(i)}\bigr)=\ell_i^{(n)}\in\bur(\per_n)$.
\end{lemma}
\begin{proof}
We pass to the representation ring. Here, since $\btr$ is a morphism of $\lambda$-rings that commutes with the induction and restriction maps, the image of the left hand side under $\btr$ is
\begin{align*}
\ind_{\per_i\times\per_{n-i}}^{\per_n}\circ\res_{\per_i\times\per_{n-i}}^{\per_i}\circ\lambda^i\bigl(\bbQ\bigl[\{1\dotsk i\}\bigr]\bigr)&&\in\rep_{\bbQ}(\per_n),
\end{align*}
and the image of the right hand side is $\lambda^i\bigl(\bbQ\bigl[\{1\dotsk n\}\bigr]\bigr)\in\rep_{\bbQ}(\per_n)$. Now $\ell_i^{(i)}\in\schur_i$ hence, by the preceding lemma, $\ind_{\per_i\times\per_{n-i}}^{\per_n}\circ\res_{\per_i\times\per_{n-i}}^{\per_i}\ell_i^{(i)}\in\schur_n$. Since also $\ell_i^{(n)}\in\schur_n$ and $\btr$ is injective on $\schur_n$, it suffices to prove that $$\ind_{\per_i\times\per_{n-i}}^{\per_n}\circ\res_{\per_i\times\per_{n-i}}^{\per_i}\Bigl(\lambda^i\bigl(\bbQ\bigl[\{1\dotsk i\}\bigr]\bigr)\Bigr)=\lambda^i\bigl(\bbQ\bigl[\{1\dotsk n\}\bigr]\bigr)\in\rep_{\bbQ}(\per_n),$$
i.e., we have to find a $\per_n$-equivariant isomorphism of $\bbQ$-vector spaces
\begin{equation*}
\varphi\colon\bbQ[\per_n]\otimes_{\bbQ[\per_i\times\per_{n-i}]}\ext^i\bbQ\bigl[\{1\dotsk i\}\bigr]\rightarrow\ext^i\bbQ\bigl[\{1\dotsk n\}\bigr].
\end{equation*}
This is straightforward. (It is done explicitly in \cite{rokaeus07}, Proposition 2.4.2.)
\end{proof}

We are now ready to prove the main theorem of this section, the formula for $\ell_i^{(n)}$. For this we first introduce a concept of degree on basis elements of $\schur_n$. Fix an $n$ and the basis $\{[\mathcal{P}^{(n)}_\mu]\}_{\mu\vdash n}$ of $\schur_n$. For $j=1\dotsk k$ where $2k<n$, we say that an element $[\ccP^{(n)}_\mu]$ in the basis is of degree $j$ if it is equal to $[\ccP^{(n)}_\nu]$ for some $\nu\vdash j$. Equivalently, this means that $n-j$ is an entry of $\mu$. Let the degree of $[\ccP^{(n)}_n]=1$ be zero and let the degree of the remaining elements of the basis be $k+1$. Because $n-j>k$ the degree is well-defined.
\begin{lemma}\label{9}
Let $[\ccP^{(n)}_\alpha]$ and $[\ccP^{(n)}_\beta]$ be of degree $m$ and $m'$ respectively, where $m+m'\leq n/2$. Then \begin{equation*}\label{5}
[\cc{P}_{\alpha}^{(n)}]\cdot[\cc{P}^{(n)}_{\beta}]=[\cc{P}^{(n)}_{\alpha,\beta}]+\text{ terms of degree $<m+m'$.}
\end{equation*}
\end{lemma}
\begin{proof}
This is a refinement of Proposition \ref{7}. Let $s=(S_1\dotsk S_l)\in\ccP^{(n)}_\alpha$ and $t=(T_1\dotsk T_{l'})\in\ccP^{(n)}_\beta$, where $S_l$ and $T_{l'}$ have $n-m$ and $n-m'$ elements respectively. Then the stabilizer of $(s,t)\in\ccP^{(n)}_\alpha\times\ccP^{(n)}_\beta$ equals the stabilizer of $(S_i\cap T_j)_{i,j}$. Let $m_{ij}=\abs{S_i\cap T_j}$ and let $\gamma$ be the tuple consisting of the $m_{ij}$. Then the transitive component of $(s,t)$ is $\ccP_\gamma^{(n)}$. Since $m_{ll'}\geq n-m-m'\geq n/2$ it follows that the degree of $[\ccP^{(n)}_\gamma]$ is $n-m_{ll'}$, and this is $\leq m+m'$ with equality if and only if $[\ccP^{(n)}_\gamma]=[\ccP^{(n)}_{\alpha,\beta}]$. 
\end{proof}
\begin{theorem}\label{80}
Let $i$ be a positive integer. Then for any $n\geq i$,
\begin{equation*}
\lambda^i\bigl(\{1\dotsk n\}\bigr)=(-1)^i\sum_{\mu\vdash i}(-1)^{\ell(\mu)}\tbinom{\ell(\mu)}{\alpha(\mu)}    \bigl[\mathcal{P}_{\mu}^{(n)}\bigr]    \in\bur(\per_n).
\end{equation*}
\end{theorem}
\begin{proof}
For $i=1$ the formula becomes $\ell_1^{(n)}=\bigl[\cc{P}_{1}^{(n)}\bigr]=[\{1\dotsk n\}]$, which is true for every $n$.

Given $i$, suppose the formula is true for every pair $(i',n)$ where $i'<i$ and $n$ is an arbitrary integer greater than or equal to $i'$. We want to show that it holds for $(i,n)$ where $n$ is an arbitrary integer greater than or equal to $i$.

Since $\ell_i^{(i)}\in\schur_i$, we have $\ell_i^{(i)}=\sum_{\mu\vdash i}a_\mu\bigl[\mathcal{P}_{\mu}^{(i)}\bigr]$,
where the $a_\mu$ are uniquely determined. Since the induction and restriction maps are additive it follows from lemma \ref{74} and lemma \ref{73} that
\begin{equation}\label{69}
\ell_i^{(n)}=\sum_{\mu\vdash i}a_\mu\bigl[\mathcal{P}_{\mu,n-i}^{(n)}\bigr]\in\schur_n
\end{equation}
for every $n\geq i$. It remains to show that $a_\mu=(-1)^i(-1)^{\ell(\mu)}\tbinom{\ell(\mu)}{\alpha(\mu)}$ for every $\mu\vdash i$. For this, fix an $n$ greater than $2i$ and the basis $\{[\mathcal{P}_\mu]\}_{\mu\vdash n}$ of $\schur_n$. We now use the notion of degree introduced before this theorem. By \eqref{69}, $\ell_i^{(n)}$ is a linear combination of elements of degree $i$. On the other hand, by the definition of $\ell_i^{(n)}$ we have
\begin{equation}\label{eq:66}
-(-1)^i\ell_i^{(n)}=\sum_{j=0}^{i-1}(-1)^j\ell_j^{(n)}s_{i-j}^{(n)}.
\end{equation}
By induction and the formula for $s_j^{(n)}$ from Proposition \ref{8}, the right hand side of \eqref{eq:66} equals
\begin{equation}\label{117}
\sum_{\mu\vdash i}[\cc{P}^{(n)}_{\alpha(\mu)}]+\sum_{j=1}^{i-1}(-1)^j\Biggl( (-1)^j    \sum_{\mu\vdash j}(-1)^{\ell(\mu)}\tbinom{\ell(\mu)}{\alpha(\mu)}\bigl[\mathcal{P}_{\mu}^{(n)}\bigr]       \Biggr)\cdot\Biggl(\sum_{\mu\vdash i-j}[\cc{P}^{(n)}_{\alpha(\mu)}]\Biggr)
\end{equation}
Since we already know that $\ell_i^{(n)}$ is zero in every degree different from $i$ it remains to compute the degree $i$ part of \eqref{117}. In this expression, for every $j$ such that $0<j<i$ we have a product of two sums, one consisting of elements of degree $j$ and the other one consisting of elements of degree less than or equal to $i-j$, for if $\mu\vdash i-j$ then $[\cc{P}_{\alpha(\mu)}]$ has degree $\leq i-j$ with equality if and only if $\mu=(1,1\dotsk 1)$, in which case $\alpha(\mu)=(i-j)$. If $[\cc{P}_\mu^{(n)}]$ has degree $j$ and $[\cc{P}^{(n)}_{\alpha(\mu')}]$ has degree $m\leq i-j$ then, by Lemma \ref{9},
$$[\cc{P}_{\mu}^{(n)}]\cdot[\cc{P}^{(n)}_{\alpha(\mu')}]=[\cc{P}^{(n)}_{\mu,\alpha(\mu')}]+\text{ terms of degree $<j+m$.}$$
Hence only the degree $i-j$ part of $\sum_{\mu\vdash i-j}[\cc{P}^{(n)}_{\alpha(\mu)}]$, i.e., $[\cc{P}^{(n)}_{i-j}]$, contributes to the degree $i$ part of \eqref{117}, which therefore equals
\begin{equation}\label{118}
[\cc{P}_i^{(n)}]+\sum_{j=1}^{i-1}\sum_{\mu\vdash j}(-1)^{\ell(\mu)}\tbinom{\ell(\mu)}{\alpha(\mu)}\bigl[\mathcal{P}_{\mu,i-j}^{(n)}\bigr].
\end{equation}
We write this as a linear combination of elements in $\{[\ccP_\nu^{(n)}]\}_{\nu\vdash i}$. Fix $\nu\vdash i$ with $\ell:=\ell(\nu)$ and $\alpha:=\alpha(\nu)=(\alpha_1\dotsk\alpha_t)$. If $\ell=1$ then $[\ccP_\nu^{(n)}]=[\ccP_i^{(n)}]$, so $[\ccP_\nu^{(n)}]$ occurs once in \eqref{118}. If $\ell>1$ then $[\ccP_\nu^{(n)}]$ occurs first when $i-j$ equals $\nu_1=\dots=\nu_{\alpha_1}$; the length of $\mu$ is then $\ell-1$ and $\alpha(\mu)=(\alpha_1-1,\alpha_1\dotsk\alpha_t)$, so the coefficient in front of $[\ccP_{\mu,i-j}^{(n)}]$ is $-(-1)^\ell\alpha_1\cdot\tfrac{(\ell-1)!}{\alpha!}$. Also, $[\ccP_\nu^{(n)}]$ occurs in \eqref{118} when $i-j$ equals $\nu_{\alpha_1+1}=\dots=\nu_{\alpha_1+\alpha_2}$, with coefficient $-(-1)^\ell\alpha_2\cdot\tfrac{(\ell-1)!}{\alpha!}$, and so on. Summing up, the coefficient in front of $[\ccP_\nu^{(n)}]$ is $-(-1)^\ell(\alpha_1+\dots+\alpha_t)\tfrac{(\ell-1)!}{\alpha!}=-(-1)^\ell\tbinom{\ell}{\alpha}$; hence \eqref{118} equals
\begin{equation*}
-\sum_{\nu\vdash i}(-1)^{\ell(\nu)}\tbinom{\ell(\nu)}{\alpha(\nu)}\bigl[\mathcal{P}_{\nu}^{(n)}\bigr].
\end{equation*}
Therefore, by \eqref{eq:66}, $\ell_i^{(n)}$ has the desired form and by induction we are through.
\end{proof}
\begin{remark}
This proof starts with noting that, as a consequence of the preceding lemmas, given $i$ it suffices to compute $\ell_i^{(n)}$ for some $n$ in order to get the formula for every $n$. We then compute $\ell_i^{(n)}$ for $n$ sufficiently large. Instead we could have computed $\ell_i^{(i)}$ by using that $\btr(\ell_i^{(i)})=[\sgn]\in\rep_{\bbQ}(\per_i)$, where $\sgn$ is the signature representation. The needed expression of $[\sgn]$ as a linear combination of permutation representations is a classical formula in the theory of representations of $\per_i$. We chose to give the above proof since it is purely combinatorial in nature.
\end{remark}

\bibliography{bibNoteOnTheBurnsideRing}
\bibliographystyle{amsalpha}

\end{document}